\documentclass[11pt]{article} 
\usepackage{amsmath,amsthm,amscd,amssymb}
\usepackage{latexsym}
\usepackage{epsf}
\usepackage{epsfig}

\newcommand{\diff}{\operatorname{Diff}}

\newcommand{\R}{\mathbb{R}}

\newcommand{\N}{\mathbb{N}}

\newcommand{\w}{\omega}

\newcommand{\ep}{\varepsilon}

\newcommand{\al}{\alpha}
\newcommand{\be}{\beta}

\newcommand{\LM}{L(M)}
\newcommand{\LN}{L(N)}
\newcommand{\Mint}{[0,\LM]}
\newcommand{\Nint}{[0,\LN]}

\newcommand{\SA}{\mathcal{A}}

\newcommand{\SB}{\mathcal{B}}
\newcommand{\SC}{\mathcal{C}}

\newcommand{\SQ}{\mathcal{Q}}

\newcommand{\pd}{\partial}

\theoremstyle{plain}
\newtheorem{theorem}{Theorem}[section]
\newtheorem{lemma}[theorem]{Lemma}
\newtheorem{corollary}[theorem]{Corollary}
\newtheorem{proposition}[theorem]{Proposition}

\theoremstyle{definition}
\newtheorem{definition}[theorem]{Definition}

\newtheorem{example}[theorem]{Example}
\theoremstyle{remark}
\newtheorem{remark}[theorem]{Remark}

\newtheorem{case[theorem]}{Case}

\title{Deformation Minimal Bending of Compact\\ Manifolds: Case of Simple Closed
Curves \footnote{Supported by NSF Grant DMS 0604331}}
\author{Oksana Bihun \footnote{Corresponding author e-mail: oksana@math.missouri.edu}
and Carmen Chicone\footnote{e-mail: carmen@math.missouri.edu}\\
University of Missouri, Columbia, USA}

\begin{document}
\maketitle
\begin{abstract}
The problem of minimal distortion bending of smooth 
 compact embedded connected Riemannian $n$-manifolds
$M$ and $N$ without boundary is made precise by defining a deformation energy
functional $\Phi$ on the set of diffeomorphisms $\diff(M,N)$. We derive the  
Euler-Lagrange equation for $\Phi$ and
determine smooth minimizers of $\Phi$ in case 
$M$ and $N$ are simple closed curves. 
\end{abstract}

\noindent {\bf MSC 2000 Classification:} 58E99\\
{\bf Key words:} minimal deformation,  
distortion minimal, geometric\\ optimization

\section{Introduction}
Two diffeomorphic compact embedded hypersurfaces admit infinitely many 
diffeomorphisms, which we view as  prescriptions for bending one hypersurface 
into the other. We ask which diffeomorphic bendings have minimal distortion 
with respect to a natural bending energy functional that will be precisely 
defined. 
We determine the Euler-Lagrange equation for the general case of hypersurfaces 
in Euclidean spaces and solve the problem for one-dimensional manifolds 
embedded in the plane. The existence of minima for the general case 
is a difficult open problem. 
An equivalent problem for a functional that measures the total energy of 
deformation due to stretching was solved in \cite{BC}. Some related 
discussions on the 
minimization problem are presented in \cite{LYS,Terz,W, YPM}.

\section{Minimal Distortion Diffeomorphisms}
\label{MNprop}
Let $M$ and $N$ denote compact, 
connected and oriented $n$-manifolds without boundary that are embedded 
in $\R^{n+1}$ and equip them with the natural Riemannian metrics  
$g_M$ and $g_N$ inherited from the usual metric of $\R^{n+1}$. These Riemannian manifolds $(M,g_M)$ and $(N,g_N)$ have the volume forms
$\omega_M$ and $\omega_N$  induced by their Riemannian metrics.
We assume that $M$ and $N$ are diffeomorphic,  denote the class of 
($C^\infty$) diffeomorphisms from $M$ to $N$ by $\diff(M,N)$, the 
(total space of the) tangent bundle of $M$ by $TM$, the cotangent bundle  by $TM^\ast$, and the sections of an arbitrary vector bundle $V$  by $\Gamma(V)$. For $h\in \diff(M,N)$, 
we use the standard notation $h^*$ for the pull-back map associated with $h$ 
and $h_*$ for its push-forward map. 

\begin{definition}
The strain tensor $S \in \Gamma(TM^\ast \otimes TM^\ast )$ corresponding to 
$h \in \diff(M,N)$ is defined to be
\begin{equation}
S=h^\ast g_N- g_M
\end{equation}
(cf. \cite{Terz}, \cite{MH}).
\end{definition}

Recall the natural bijection between covectors in 
$T^*M$ and vectors  in $TM$ (see \cite{KN}): To each covector  
$\alpha_p \in T_pM^\ast$ assign the 
vector  $\alpha^\#_p \in T_p M$ that is implicitly defined by the relation
$$
\alpha_p=(g_M)_p(\alpha^\#_p,\cdot).
$$
Using this correspondence, we introduce the Riemannian metric $g_M^\ast$ on   $TM^\ast$ by  
$$
g_M^\ast(\al, \be)=g_M(\al^\#, \be^\#),
$$
where the base points are suppressed. 

There is a natural Riemannian metric $G$ on $TM^\ast \otimes TM^\ast$ given by $G=g^\ast_M \otimes
g^\ast_M$. To compute this metric in local coordinates, 
let $(U, \phi)$ be a local coordinate system on $M$.  Using the coordinates of $\R^n$,  the map $\phi:U\to \R^n$ can be expressed
in the form
$$
\phi(p)=\big(x^1(p), \ldots, x^n(p)\big).
$$
As usual $\big(x^1(p), \ldots, x^n(p)\big)$ are the local coordinates of $p \in M$
and the $n$-tuple of functions $(x^1, x^2, \ldots, x^n)$ is the local coordinate
system with respect to $(U, \phi)$. Because $\phi$ is a homeomorphism from $U$
onto $\phi(U)$, we identify  $p \in U$ and $\phi(p) \in \R^n$ via $\phi$.
Let us define $\Big(\frac{\pd}{\pd
x^i}\Big)_p=\frac{\pd \phi^{-1}}{\pd x^i}\big(\phi(p)\big)$.
The set of
vectors $\Big((\frac{\pd}{\pd x^1})_p, \ldots, (\frac{\pd}{\pd x^n})_p\Big)$ 
forms a basis of the tangent space $T_pM$. Its dual basis $\big((dx^1)_p, \ldots,
(dx^n)_p\big)$ is a basis of $T_pM^\ast$, i.e.
$$
(dx^i)_p\Big(\big(\frac{\pd}{\pd
x^j}\big)_p \Big)=\delta^i_j, \,\,\,1 \leq i,j \leq n.
$$
Using the Einstein summation convention, a tensor 
$B \in \Gamma(TM^\ast \otimes TM^\ast)$ has local coordinate
representation
$B=b_{ij}dx^i\otimes dx^j$,
where $b_{ij}=B(\pd/\pd x^i,\pd/\pd x^j)$. 
 The local coordinate representation 
of the Riemannian metric $G$ is  
\begin{eqnarray}
\nonumber
G(B,B)&=&b_{ij}b_{kl}g^\ast_M(dx^i,dx^k)g^\ast_M(dx^j,dx^l)\\
\label{Gformula}
&=&
b_{ij}b_{kl}[g_M]^{ik} [g_M]^{jl},
\end{eqnarray}
where $[g_M]^{ij}$ is the $(i,j)$ entry of the inverse matrix of 
$\big([g_M]_{ij}\big)$.

\begin{definition}
\label{DefPhi}
The deformation energy functional $\Phi:\diff(M,N)\to \R_+$ is defined to be
\begin{equation}
\label{phidef1}
\Phi(h)=\int_M G(h^\ast g_N- g_M, h^\ast g_N- g_M)\,\w_M.
\end{equation}
\end{definition}

The following invariance property of the functional
 $\Phi$ is obvious because the isometries of $\R^{n+1}$ are compositions
of translations and rotations, which produce no deformations.

\begin{lemma}
\label{l1}
If $k \in \diff(N)$ is an isometry of $N$ (i.e. $k^\ast g_N=g_N$), then
$\Phi(k\circ h)=\Phi(h)$.
\end{lemma}

\section{The First Variation}

We will compute the Euler-Lagrange equation for the deformation energy
functional $\Phi$. To do this,  we will consider smooth variations. 

\begin{definition}
A $C^\infty$ function $F(t,p)=h_t(p)$ defined on $(-\ep, \ep) \times
M$ is called a smooth variation of a diffeomorphism $h \in \diff(M,N)$ if 
\begin{enumerate}
\item $h_t \in \diff(M,N)$ for all $t \in (-\ep, \ep)$ and
\item $h_0=h$.
\end{enumerate}
\end{definition}

The tangent space $T_h \diff(M,N)$ is identified with the set 
$\Gamma(h^{-1}TN)$ of all the
smooth sections of the induced bundle $h^{-1}TN$ with fiber $T_{h(p)}N$ 
over the point $p$ of the manifold $M$ (cf. \cite{Nishikawa}).
Indeed, each smooth variation $F:(-\ep, \ep) \times
M \to N$ corresponds to  a curve $t \mapsto F(t,p)=h_t(p)$ in $\diff(M,N)$.

\begin{definition}
Let $F:(-\ep, \ep) \times M \to N$ be a smooth variation of a diffeomorphism $h
\in \diff(M,N)$. The variational vector field $V \in \Gamma(h^{-1}TN)$ is
defined by 
$$
V(p)=\frac{d}{dt}h_t\big|_{t=0}(p)=\frac{\pd}{\pd t}F(0,p)
$$
for  $p \in M$.
\end{definition}

Since the tangent space $T_h \diff(M,N)$ consists of all the variational vector fields of the diffeomorphism $h$, it follows that $T_h
\diff(M,N)$ is a subset of $\Gamma(h^{-1}TN)$. On the other hand, suppose that a
 vector field $V \in \Gamma(h^{-1}TN)$ is given. We can easily construct a
variation of $h$ with the variational vector field $V$. Indeed, let $\psi_t$ be
the flow of the vector field $X=V \circ h^{-1} \in \Gamma(TN)$. The smooth
variation $F(t,p)= \psi_t \circ h(p)$ of the diffeomorphism $h \in \diff(M,N)$
has the variational vector field $V(p)=\frac{d}{dt}( \psi_t \circ h)(p)= X \circ
h (p)= V(p)$ as required. Hence, \[T_h \diff(M,N)=\Gamma(h^{-1}TN).\]

We will consider all variations of $h \in \diff(M,N)$ of the
form $F(t,p)=h \circ \phi_t (p)$, where $\phi_t$ is the flow of a vector field
$X \in \Gamma(TM)$. The variational vector field corresponding to the variation $F$ is
$V=h_{\ast} X$. Since $h$ is a diffeomorphism, it is easy to see that the
variational vector fields of the variations of the form $h \circ \psi_t$ exhaust all
possible variational vector fields. 

Let us restrict the domain of the functional $\Phi$ to 
$\diff(M,N)$. The diffeomorphism $h$ is a critical point of $\Phi$ if 
\begin{equation}
\label{EL3}
\frac{d}{dt} \Phi(h \circ \phi_t)|_{t=0 }=
D\Phi(h) h_\ast Y=\int_M G(h^\ast g_N-g_M, L_Y h^\ast g_N)=0
\end{equation}
for all $Y \in \Gamma(TM)$, where $L_Y$ denotes the Lie derivative in the
direction $Y$.

Let $\be \in \Gamma(TM^\ast \otimes TM^\ast)$ have the
  local representation
$\be_{ij} dx^i \otimes dx^j$.
We will use the following formula for the 
components of the Lie derivative $L_X\be$ of $\be$ in the direction of the vector field $X$:
\begin{equation}
\label{LXgNij}
[L_X \be]_{ij}=X^k\frac{\pd \be_{ij}}{\pd x^k}+
\be_{kj}\frac{\pd X^k}{\pd x^i}+\be_{ik}\frac{\pd X^k}{\pd x^j}.
\end{equation}

\section{Solution for One Dimensional Manifolds}
\label{1dim}
In this section  
$M$ and $N$ are smooth simple closed curves in $\R^2$. Their arclengths are denoted $\LM$ and $\LN$ respectively, and they are supposed to have base points $p\in M$ and $q\in N$. 
We will determine the minimum of the functional 
\begin{equation}
\label{e21}
\Phi(h)=\int_M G\big(h^\ast g_N-g_M, h^\ast g_N-g_M\big)\w_M
\end{equation}
over the admissible set 
\begin{equation}
\label{e211}
\SA=\{h \in \diff(M,N):h(p)=q\}.
\end{equation} 

There exist unique arc length parametrizations $\gamma: \Mint\to M$ and
$\xi:\Nint \to N$ of $M$ and $N$ respectively, which correspond to the positive orientations of the curves $M$ and
 $N$ in the plane, and are such that $\gamma(0)=p$, 
$\xi(0)=q$.
Notice that $[g_M]_{11}(t)=|\dot{\gamma}(t)|^2=1=[g_M]^{11}(t)$ for $t \in [0,L(M)]$
and $[h^\ast g_N]_{11}(t)=|Dh\big(\gamma(t)\big)\dot{\gamma}(t)|^2$. Using
formula (\ref{Gformula}) for the metric $G$, we can rewrite functional 
(\ref{e21}) in local coordinates:
\begin{equation}
\label{e22}
\Phi(h)=\int_0^{\LM} 
\Big(\big|Dh\big(\gamma(t)\big)\dot{\gamma}(t) \big|^2-1\Big)^2\w_M.
\end{equation}
Let us denote the local representation of a diffeomorphism $h \in \diff(M,N)$ by
$u=\xi^{-1}\circ h \circ \gamma$. The function $u$ is a diffeomorphism on the open
interval $\big(0,\LM\big)$ and can be continuously extended to the closed interval
$[0,\LM]$ as follows. If $h$ is orientation preserving, we can extend $u$ to a
continuous function on $\Mint$ by defining $u(0)=0$ and $u(\LM)=\LN$. In this case
$\dot{u}>0$. If $h$ is orientation reversing, we define 
$u(0)=\LN$ and $u(\LM)=0$.

Since
$$
\big|\frac{d}{dt}(h\circ \gamma)(t)\big|^2=
\big|\frac{d}{dt}(\xi\circ u)(t)\big|^2=\dot{u}^2(t)
\big|\dot{\xi}\big(u(t)\big)\big|^2=\dot{u}^2(t)
$$
for $t \in \big(0,\LM\big)$,
the original problem of the minimization of functional (\ref{e21}) can be
reduced to the minimization of the functional
\begin{equation}
\label{e23}
\Psi(u)=\int_0^{\LM} (\dot{u}^2-1)^2dt
\end{equation}
over the admissible sets
$$\SB=\Big\{u\in C^2\big(\Mint,\Nint\big): u(0)=0, u(\LM)=\LN \Big\}$$ 
and 
$$\SC=\Big\{u\in C^2\big(\Mint,\Nint\big): u(0)=\LN, u(\LM)=0 \Big\}.$$
The minima will be shown to correspond to diffeomorphisms in $\diff(M,N)$. 

\begin{lemma}
\label{ll1}
Suppose that  $\LN \geq \LM$.
\begin{itemize}
\item[(i)]
The function $v(t)=\LN/\LM t$, where $t \in \Mint$, is the unique minimum of the
functional $\Psi$ over the admissible set $\SB$.
\item[(ii)] The function $w(t)=-\LN/\LM t+\LN$ , where $t \in \Mint$, is the unique minimum of the
functional $\Psi$ over the admissible set $\SC$.
\end{itemize}
\end{lemma}
\begin{proof}
Since the proofs of (i) and (ii) are almost identical, we will only present the
proof of the statement (i).

The Euler-Lagrange equation for functional (\ref{e23}) is
\begin{equation}
\label{e215}
4\ddot{u}(3\dot{u}^2-1)=0.
\end{equation}
The only solution of the above equation that belongs to the admissible set $\SB$
is $v(t)=\frac{\LN}{\LM}t$, where $t \in \Mint$. Note that $v$ 
corresponds to a diffeomorphism in $\diff(M,N)$. 

We will show that the critical
point $v$ minimizes the functional $\Psi$; that is, 
\begin{equation}
\label{a5455}
\Psi(u)\geq\Psi(v)=\frac{(\LN^2-\LM^2)^2}{\LM^3}
\end{equation}
for all $u \in \SB$. Using H\"older's inequality
$$
\LN=u(\LM)=\int_0^{\LM} \dot{u}(s)\,ds\leq 
\big[\LM\int_0^{\LM} \dot{u}^2(s)\,ds\big]^{1/2},
$$
we have that
$$
\frac{\LN^2}{\LM}\leq \int_0^{\LM} \dot{u}^2(s)\,ds.
$$
Thus, in view of the hypothesis that  $\LN \geq \LM$, 
\begin{eqnarray}
\nonumber
\int_0^{\LM}(\dot{u}^2(s)-1)\,ds&=&\int_0^{\LM} \dot{u}^2(s)\,ds-\LM \\
\label{e24}
&\geq& \frac{\LN^2-\LM^2}{\LM}\geq 0.
\end{eqnarray}
After squaring both sides of inequality (\ref{e24}), we obtain the inequality
\begin{equation}
\label{e25}
\big(\int_0^{\LM}(\dot{u}^2(s)-1)\,ds\big)^2 \geq \frac{(\LN^2-\LM^2)^2}{\LM^2}.
\end{equation}
Applying H\"older's inequality to $\Phi(u)$ and taking into account inequality 
\eqref{e25}, we obtain inequality \eqref{a5455}.
Hence, the function $v(t)=\LN/\LM t$, where $t \in \Mint$, minimizes the
functional $\Psi$ over the admissible set $\SB$.
\end{proof}

\begin{remark}
\label{remEL}
Let us write the Euler-Lagrange equation (\ref{EL3}) for the one-dimensional
case and compare it with equation (\ref{e215}).

Recall that
$$
[g_M]_{11}(t)=1, \qquad [h^\ast g_N]_{11}(t)=\dot{u}(t)^2
$$
and use formula \eqref{LXgNij} to compute  
$$ [L_Y h^\ast g_N]_{11}(t)=2\dot{u}(t)\big(\ddot{u}(t)y(t)+
\dot{y}(t)\dot{u}(t)\big)=
2\dot{u}(t)\frac{d}{dt}\big(\dot{u}(t)y(t)\big),
$$
where $y(t)$ is the local coordinate of the vector field $Y=y \frac{\pd}{\pd t}$, 
i.e., $y$ is a smooth periodic function
on $\Mint$, which can be taken to be in $C_c^{\infty}(\Mint)$.
Using the pervious computation and formulas (\ref{Gformula}) and (\ref{EL3}), we obtain the
following Euler-Lagrange equation:
$$
\int_0^{\LM} (\dot{u}^2-1)\dot{u}\frac{d}{dt}(\dot{u}y) \,dt=-
\int_0^{\LM} \frac{d}{dt}\big((\dot{u}^2-1)\dot{u}\big)\dot{u}y \,dt=0
$$
for all $y \in C_c^{\infty}(\Mint)$. The latter equation yields
\begin{equation}
\label{EL1dim}
\frac{d}{dt}\big((\dot{u}^2-1)\dot{u}\big)\dot{u}=\dot{u} 
\ddot{u}(3 \dot{u}^2-1)= 0,
\end{equation}
which has the same solutions in the admissible sets $\SB$ and $\SC$ as equation 
(\ref{e215}).
\end{remark}

\begin{proposition}
\label{cor1dim}
Suppose that $M$ and $N$ are smooth simple closed curves in $\R^2$ with arc lengths $\LM$ and $\LN$ and base points $p\in M$ and $q \in N$; 
 $\gamma$ and $\xi$ are  arc length parametrizations of $M$ and $N$ with $\gamma(0)=p$ and $\xi(0)=q$ that induce positive orientations; and, the functions $v$ and $w$ are as in
lemma \ref{ll1}.
If $\LN \geq\LM$, then
the functional $\Phi(h)$ defined in display (\ref{e21}) has exactly two
minimizers in the admissible set $$\SA=\{h \in \diff(M,N):h(p)=q\}:$$
the orientation preserving minimizer 
\[
h_1=\xi \circ v \circ \gamma^{-1}
\]
and the orientation reversing minimizer
\[
h_2=\xi \circ w \circ \gamma^{-1}
\]
(where we consider $\gamma$ as a function defined on  
$\big[0,\LM \big)$ 
so that $\gamma^{-1}(p)=0$). Moreover, 
the minimal value of the functional $\Phi$ is
\begin{equation}
\Phi_{min}=\frac{(\LN^2-\LM^2)^2}{\LM^3}.
\end{equation}
\end{proposition}

\begin{example}
For $R>0$, the radial map $h:\R^2 \to \R^2$ is defined to be $h(z)=Rz$. 
If $M$ is a simple closed curve, $N:=h(M)$ and $R>1$, then $h$ is a minimum 
of $\Phi$ on $\diff(M,N)$.
To see this fact,
let $\gamma(t)=\big(x(t), y(t)\big)$, $t \in \Mint$, be an arc length parametrization of $M$. It
is easy to see that $\xi(t)=R \big(x(t/R), y(t/R)\big)$, 
$t \in [0, R\LM]$ parametrizes $N=h(M)$ by its arc length. By proposition
\ref{cor1dim}, the
minimizer $h_1$  is 
\begin{eqnarray*}
h_1(z)=\xi\big(v \circ
\gamma^{-1}(z)\big)&=&\xi\big(R \gamma^{-1}(z)\big)\\
&=&\xi(Rt)=R \gamma(t)=Rz
\end{eqnarray*}
for all $z \in M$. Hence, $h_1$ is the radial map. 
\end{example}

\begin{lemma}
\label{le:nomin}
If $\LN<\LM$,  then the functional $\Psi$ has no minimum in the 
admissible set $\SB$.
\end{lemma}
\begin{proof}
Let $\phi:[0,\LM] \to \R$ be a continuous piecewise linear function 
such that $\phi(0)=0$, $\phi(\LM)=\LN$, and  $\dot{\phi}(t)=\pm 1$ whenever 
$t\in (0,\LM)$ and the derivative is defined. The graph of $\phi$ looks
like a zig-zag.
It is easy to see that $\phi$ is an element of the Sobolev space $W^{1,4}(0,\LM)$ 
(one weak derivative in the Lebesgue space $L^4$). By the standard 
properties of  $W^{1,4}(0,\LM)$ with its usual
norm $\|\cdot\|_{1,4}$, there exists a
sequence of smooth functions $\phi_k \in C^\infty[0, \LM]$ (satisfying the boundary conditions $\phi_k(0)=0$ and $\phi_k(\LM)=\LN$) such that
$\|\phi_k-\phi\|_{1,4} \to 0$ as $k \to \infty$. Moreover, there is some 
constant $C>0$ such that 
$\int_{0}^{\LM} (\dot{\phi}_k^2-\dot{\phi}^2)^2 \,dx \leq C
\|\phi_k-\phi\|_{1,4}^2$. It is easy to see that   
$$|\Phi(\phi_k)-\Phi(\phi)|\leq C_1 \|\phi_k-\phi\|_{1,4}$$ for some constant
$C_1>0$. Taking into account the equality $\Psi(\phi)=0$, we conclude that $\Psi(\phi_k)
\to 0$ as $k \to \infty$. Thus, $\{\phi_k\}_{k=1}^\infty$ is 
a minimizing sequence for the
functional $\Psi$ in the admissible set $\SB$. On the other hand, there is no
function $f \in \SB$ such that $\Psi(f)=0=\inf_{g\in \SB} \Psi(g)$. Therefore,
if $\LN<\LM$, the functional $\Phi$ has no minimum in the admissible set $\SB$.
\end{proof}

\begin{corollary}
\label{cor1}
If $\LN<\LM$, then the functional $\Phi$ has no minimum in the admissible set
$$
\SQ=\{h \in C^2(M,N): h \,\,\,is\,\,\, orientation\,\,\, preserving
\,\,\, and\,\,\, h(p)=q \}.
$$
\end{corollary}

Let us interpret the result of Lemma~\ref{le:nomin}. Let 
$h=\xi\circ\phi\circ \gamma^{-1}$, where $\phi:[0,\LM]\to \R$ is defined in
the proof of Lemma~\ref{le:nomin} and
$\gamma$, $\xi$ are arc length (positive orientation)  
parametrizations of the curves $M$ and $N$ viewed as periodic functions on $\R$.
In case $\LN<\LM$, the action of the function $h$  on 
the curve $M$ can be described as follows.  
The curve $M$ is cut into 
segments $\{M_i\}_{i=1}^k$, $k \in \N$, such that $\dot{\phi}$ has a constant value ($1$ or $(-1)$) on 
$\gamma^{-1}(M_i)$. Each segment $M_i$ is  wrapped 
around the curve $N$ counterclockwise or clockwise depending on whether
$\dot{\phi}$ equals $1$ or $(-1)$ on $\gamma^{-1}(M_i)$ respectively.
Since $\LN$ is less than $\LM$, some 
points of $N$ will be covered by segments of $M$ several times. During
this process, the segments of the curve $M$ need not be stretched. Hence,
as  measured by the functional $\Phi$, no strain is produced, i.e. $\Phi(h)=0$.

The statement of corollary \ref{cor1} leaves open an interesting question:  
Does the functional $\Phi$  have a minimum  in the admissible set $\SA$?
 Some results in this direction are presented in the next section.

\section{Second variation}

We will derive a necessary condition for a
diffeomorphism $h \in \diff(M,N)$ to be a minimum of the functional $\Phi$.
Let $h_t=h \circ \phi_t$ be a family of diffeomorphisms in $\diff(M,N)$, where
$\phi_t$ is the flow of a vector field $Y \in \Gamma(TM)$.
Using the Lie derivative formula (see \cite{AMR}), we derive the equations
$\frac{d}{dt}(h_t^\ast g_N)=\phi_t^\ast L_Y h^\ast g_N$ and 
$\frac{d}{dt}(\phi_t^\ast L_Y h^\ast g_N)=\phi_t^\ast L_Y L_Y h^\ast g_N$. 
 If there exists 
$\delta>0$ such that $\Phi(h_t)> \Phi(h)$ for all $|t|<\delta$ and for all 
variations $h_t$ of $h$, then $h$ is called a relative minimum of $h$.
If $h \in \diff(M,N)$ is a relative minimum of $\Phi$,
then $\frac{d^2}{dt^2}\Phi(h_t)|_{t=0}>0$.

Using the previous computations of Lie derivatives, the 
second variation of $\Phi$ is
\begin{eqnarray}
\label{2var}
\frac{1}{2} \frac{d^2}{dt^2} \Phi(h_t)|_{t=0}&=&\int_M G(L_Y h^\ast g_N, 
L_Y h^\ast g_N) \w_M \\
\nonumber
&+&\int_M G(L_Y L_Yh^\ast g_N, h^\ast g_N-g_M)\w_M.
\end{eqnarray}

\begin{lemma}
\label{2ndVar}
Let $M$ and $N$ be simple closed curves parametrized by functions
$\gamma$ and $\xi$ satisfying all the
properties stated in lemma \ref{cor1dim}.
If $h \in \diff(M,N)$ minimizes the functional $\Phi$ in the admissible set 
$\SA$, then the local representation $u=\xi^{-1}\circ h \circ \gamma$ of 
$h$ satisfies the inequality
\begin{equation}
\label{2var1}
\dot{u}^2(t)\geq\frac{1}{3}
\end{equation}
for all $t \in \big(0,\LM\big)$.
\end{lemma}
\begin{proof}
Using formula (\ref{LXgNij}), we compute
$$
[L_Y h^\ast g_N]_{11}=2(\dot{u}\ddot{u}y+\dot{u}^2 \dot{y})
$$
and
$$
[L_Y L_Y h^\ast g_N]_{11}=2(\ddot{u}^2 y^2+ \dot{u}\, \dddot{u}\, y^2+5 \dot{u}
\,\ddot{u}\,\dot{y}\,y +\dot{u}^2 \ddot{y}\, y +2\dot{u}^2
\dot{y}^2).
$$
Substituting the latter expressions into formula (\ref{2var}), we obtain the 
necessary condition 
\begin{eqnarray*}
W:=4\int_0^{\LM}\dot{u}^4 \,\dot{y}^2 \,dt&+&4
\int_0^{\LM}\dot{u}^2 (\dot{u}^2-1)\,\dot{y}^2 \,dt\\
&+&2
\int_0^{\LM}\dot{u}^2 (\dot{u}^2-1)\,y \,\ddot{y} \,dt +\ldots \geq 0,
\end{eqnarray*}
where the integrands of the omitted terms all contain the factor $y$.
After integration by parts, we obtain the inequality
\begin{eqnarray}
\nonumber
W=\int_0^{\LM}\Big(4\dot{u}^4 &+&4
\dot{u}^2 (\dot{u}^2-1)\\
\label{2in}
&-&2\dot{u}^2 (\dot{u}^2-1)\big)
\dot{y}^2\,dt +\ldots \geq0.
\end{eqnarray}

Define $y(t)=\ep \rho\big(\frac{t}{\ep}\big)\zeta(t)$, where
$\rho(t)$ is a periodic ``zig-zag'' function defined by the expressions
\begin{equation}
\label{zigzag}
\rho(t)=\left\{
\begin{array}{rl}
t, &\mbox{if\ }0\leq t<1/2,\\
1-t, &\mbox{if \ } 1/2 \leq t <1,
\end{array}
\right.
\end{equation}
and $\rho(t+1)=\rho(t)$, $\zeta \in C_c^\infty\big(0,\LM\big)$. Notice that
$\dot{\rho}^2=1$ almost everywhere on $\R$ and $\dot{y}^2=\zeta^2 +O(\ep)$ when $\ep \to 0$. 
Substitute  $y$ into inequality (\ref{2in}) and pass to  the limit as $\ep \to
0$. All the omitted  terms in the expression for $W$  tend to zero
because they contain $y$ as a factor. Hence, we have the inequality
\begin{eqnarray}
\nonumber
W&=&\int_0^{\LM}\big(4\dot{u}^4 +2
\dot{u}^2 (\dot{u}^2-1)\big)
\zeta^2\,dt \geq 0,
\end{eqnarray}
which (after a standard bump function argument) reduces to the
inequality 
\begin{equation}
\dot{u}^2\geq1/3
\end{equation}
as required.
\end{proof}

\begin{proposition}
If $M$ and $N$ are simple closed curves such that their corresponding arc
lengths $\LM$ and $\LN$ satisfy the inequality 
$\frac{\LN}{\LM}<\frac{1}{\sqrt{3}}$,  then the functional $\Phi$ has 
no minimum in the admissible set $\SA$.
\end{proposition}
\begin{proof}
If $h \in \diff(M,N)$ is a minimum of the functional 
$\Phi$, then $h$ satisfies the Euler-Lagrange equation \eqref{EL3}. Let 
$\gamma$ and $\xi$ be parametrizations of the curves $M$ and $N$ with all the
properties stated in corollary \ref{cor1dim}. By remark
\ref{remEL}, the local representation  $u=\xi^{-1}\circ h \circ \gamma$ of $h$ 
satisfies the ordinary differential equation \eqref{EL1dim} on $(0,\LM)$. In
addition, $u$ must satisfy the boundary conditions $u(0)=0, u(\LM)=\LN$ or 
$u(0)=\LN, u(\LM)=0$. Hence either $u(t)=\LN/\LM t$ or $u(t)=-\LN/\LM t+\LN$. 
Since $h$ minimizes $\Phi$, by lemma \ref{2ndVar} $\dot{u}^2\geq1/3$, or,
equivalently, $\LN/\LM\geq\frac{1}{\sqrt{3}}$. This contradicts the assumption of
the theorem.
\end{proof}

\end{document}